\documentclass[a4paper,12pt]{amsart}
\usepackage{amssymb}
\usepackage{ifthen}
\usepackage[english]{babel}
\usepackage{cite}
\usepackage[dvips]{graphicx}
\nonstopmode \numberwithin{equation}{section}
\usepackage{amssymb}
\usepackage{ifthen}
\usepackage{graphicx}
\usepackage{amsmath}
\usepackage[T1]{fontenc} 
\usepackage[utf8]{inputenc}
\usepackage[usenames,dvipsnames]{color}
\usepackage{color}
\usepackage[english]{babel}
\usepackage{fancyhdr}
\usepackage{fancybox}
\usepackage{tikz}
\usepackage{bigints}

\nonstopmode \numberwithin{equation}{section}
\theoremstyle{plain}
\newtheorem{thm}[equation]{Theorem}
\newtheorem{rem}[equation]{Remark}
\newtheorem{cor}[equation]{Corollary}
\newtheorem{lem}[equation]{Lemma}

\newtheorem{prop}[equation]{Proposition}

\newtheorem{conj}{Conjecture}

\theoremstyle{definition}
\newtheorem{defn}{Definition}[section]

\newtheorem{prob}{Problem}

\newcounter{minutes}
\divide\time by 60
\newcounter{hours}
\multiply\time by 60
\addtocounter{minutes}{-\time}

\newcounter {own}
\def\theown {\thesection       .\arabic{own}}

\newenvironment{pf}[1][]{%
	\vskip 3mm
	\noindent
	\ifthenelse{\equal{#1}{}}%
	{{\slshape Proof. }}%
	{{\slshape #1.} }%
}%
{\qed\bigskip}

\newcounter{alphabet}

\newcommand{\IR}{{\mathbb R}}
\newcommand{\ID}{{\mathbb D}}

\newcommand{\real}{{\operatorname{Re}\,}}

\def\be{\begin{equation}}
	\def\ee{\end{equation}}

\newcommand{\bee}{\begin{enumerate}}
	\newcommand{\eee}{\end{enumerate}}

\newcommand{\blem}{\begin{lem}}
	\newcommand{\elem}{\end{lem}}
\newcommand{\bthm}{\begin{thm}}
	\newcommand{\ethm}{\end{thm}}
\newcommand{\bcor}{\begin{cor}}
	\newcommand{\ecor}{\end{cor}}
\newcommand{\beg}{\begin{examp}}
	\newcommand{\eeg}{\end{examp}}
\newcommand{\begs}{\begin{examples}}
	\newcommand{\eegs}{\end{examples}}

\newcommand{\bdefn}{\begin{defn}}
	\newcommand{\edefn}{\end{defn}}

\newcommand{\bprob}{\begin{prob}}
	\newcommand{\eprob}{\end{prob}}
\newcommand{\bei}{\begin{itemize}}
	\newcommand{\eei}{\end{itemize}}

\newcommand{\bcon}{\begin{conj}}
	\newcommand{\econ}{\end{conj}}
\newcommand{\bcons}{\begin{conjs}}
	\newcommand{\econs}{\end{conjs}}
\newcommand{\bprop}{\begin{prop}}
	\newcommand{\eprop}{\end{prop}}
\newcommand{\br}{\begin{rem}}
	\newcommand{\er}{\end{rem}}
\newcommand{\brs}{\begin{rems}}
	\newcommand{\ers}{\end{rems}}
\newcommand{\bo}{\begin{obser}}
	\newcommand{\eo}{\end{obser}}
\newcommand{\bos}{\begin{obsers}}
	\newcommand{\eos}{\end{obsers}}
\newcommand{\bpf}{\begin{pf}}
	\newcommand{\epf}{\end{pf}}
\newcommand{\ba}{\begin{array}}
	\newcommand{\ea}{\end{array}}
\newcommand{\beq}{\begin{eqnarray}}
	\newcommand{\beqq}{\begin{eqnarray*}}
		\newcommand{\eeq}{\end{eqnarray}}
	\newcommand{\eeqq}{\end{eqnarray*}}

\newcommand{\ds}{\displaystyle}

\begin{document}
	
	\title{Region of variability for certain subclass of univalent functions}
	
	\author{Jnana Preeti Parlapalli}
	\address{Jnana Preeti Parlapalli,
		Department of Computer Science,
		SRM University,
		Amaravati, Andhra Pradesh 522502, India.}
	\email{preeti.parlapalli@gmail.com} 
	
	\author{Vasudevarao Allu}
	\address{Vasudevarao Allu,
		School of Basic Science,
		Indian Institute of Technology Bhubaneswar,
		Bhubaneswar-752050, Odisha, India.}
	\email{avrao@iitbbs.ac.in}
	
	\subjclass[{AMS} Subject Classification:]{Primary 30C45, 30C50, 30C80}
	
	\def\thefootnote{}
	\footnotetext{ {\tiny File:~\jobname.tex,
			printed: \number\year-\number\month-\number\day,
			\thehours.\ifnum\theminutes<10{0}\fi\theminutes }
	} \makeatletter\def\thefootnote{\@arabic\c@footnote}\makeatother
	
	\begin{abstract} 
		Let $\mathbb{D}:=\{z\in \mathbb{C}: |z|<1\}$ be the unit disk. For $0<\alpha <1$, let $f_{\alpha}(z)=z/(1-z^\alpha)$ for $z \in \mathbb{D}$. We consider the class $\mathcal{F}$ of analytic functions $f_{\alpha}$ which satisfy $\real \left(1+zf''_{\alpha}(z)/f'_{\alpha}(z)\right) > \beta$ for $0<\beta<1$. In this paper, we determine the
		region of variability of $\log f'_{\alpha}(z_0)$ for fixed $z_{0} \in \mathbb{D}$ when $f$
		varies over the class ${\mathcal F}(\lambda):=\{f_{\alpha} \in \mathcal{F}: f_{\alpha}(0)=0, f'_{\alpha}(0)=1 \, \mbox{and} \, f''_{\alpha}(0)=2\lambda (1-\beta) \,\,\, \mbox{for} \,\, 0\leq \lambda \leq 1\}$.
	\end{abstract}
	
	\maketitle
	\pagestyle{myheadings}
	\markboth{Jnana Preeti Parlapalli and Vasudevarao Allu}{Region of variability for certain subclass of univalent functions}

\section{Introduction}
Let  $\mathbb{D} := \{ z: \,  |z| < 1 \}$ be the unit disk in the
complex plane $\mathbb{C}$ and $\mathcal{H}$ denote the space of
all analytic functions on $\ID$. Here we think of $\mathcal{H}$ as
a topological vector space endowed with the topology of uniform
convergence over compact subsets of $\mathbb{D}$. A function $f$ is said to be univalent if it is one-to-one. Further, let
$\mathcal{A} := \{ f \in \mathcal{H}:\, f(0) = f'(0) - 1 = 0 \}$
and $\mathcal{S}$ denote the class of {\em univalent} functions in
$\mathcal{A}$. A domain $\Omega \subseteq \mathbb{C}$ is called a starlike domain with respect to a point $z_{0} \in \Omega$ if the line segment $tw + (1-t)z_{0}$ \, joining $z_0$ to any other point $w$ in $\Omega$ lies entirely in $\Omega$. If $z_0=0$ then the domain $\Omega$ is called a starlike domain. A function $f\in\mathcal{S}$ is said to be starlike with respect to $w_{0} \in f(\mathbb{D})$ if the image domain $f(\mathbb{D})$ is starlike with respect to $w_{0}$. If $w_{0}=0$, $f$ is simply called a starlike function.
The class of univalent starlike  functions is denoted by ${\mathcal S}^*$. 
A domain $\Omega \subseteq \mathbb{C}$ is said to be convex if the line segment joining any two points of the domain $\Omega$ is contained in $\Omega$.  A function $f$ is said to be convex if the image of the domain $\Omega$ under $f$ is a convex domain. The class of univalent convex  functions is denoted by ${\mathcal C}$.\\

 A function $f \in \mathcal{A}$ is called {\em close-to-convex} if there exists a convex (univalent) function $g$ and a number $\phi \in \IR$ such that ${\rm Re\,}(e^{i\phi }f'(z)/g'(z))>0$
for $z\in \ID$. Each univalent starlike function $f$ is characterized by the
analytic condition ${\rm Re\,}(zf'(z)/f(z))>0$ in $\ID$. Due to Alexender's theorem 
it is known that $zf'$ is starlike if, and only if, $f$ is convex (see \cite{duren-1983}). In $1952$, Kaplan \cite{Kaplan-1952} proved that
every close-to-convex function is univalent in $\ID$. A function $f \in \mathcal{S}$ is said to be starlike function of order $\alpha, -1/2 \leq \alpha < 1$ if it satisfies the condition
\begin{equation*}
	{\rm Re\,}\left(\frac{z f'(z)}{f(z)}\right)> \alpha \,\,\,\,\, \mbox{for} \quad z\in\mathbb{D}.
\end{equation*}
Let $\mathcal{S}^{*}(\alpha)$ be the class of all starlike functions of order $\alpha$ in $\mathbb{D}$. A function $f \in \mathcal{S}$ is said to be convex function of order $\alpha, -1/2 \leq \alpha < 1$ if it satisfies the condition
\begin{equation*}
	{\rm Re\,}\left(1+\frac{z f''(z)}{f'(z)}\right)> \alpha\,\,\,\, \mbox{for} \quad z\in\mathbb{D}.
\end{equation*}
Let $\mathcal{C}(\alpha)$ be the class of all convex functions of order $\alpha$ in $\mathbb{D}$.

\vspace{2mm}
Let $\mathcal{F}\subset \mathcal{A}$ and $z_0 \in \mathbb{D}$. Then upper and lower estimates of the form 
$$
K_{1} \leq f(z_0)\leq K_{2},\,\,\, M_1\leq |f'(z_0)|\leq M_2,\,\,\, m_1\leq \mbox{Arg}f'(z_0)\leq m_2 \,\,\mbox{for all} \,\, f\in \mathcal{F},
 $$
are respectively called a growth theorem, a distortion theorem and a rotation theorem at $z_0$ for $\mathcal{F}$, where $K_{i}, M_i$ and $m_i$ $(i=1,2)$ are some non-negative constants. These estimates deal only with the absolute values of $f(z_0)$ and $f'(z_0)$ or with the argument of $f'(z_0)$. If one wants to study the exact value of $f(z_0)$ and $f'(z_0)$, then it is necessary to consider the region of variability of $f(z_0)$ and $f'(z_0)$ when $f$ ranges over the class $\mathcal{F}$.
\vspace{2mm}

Region of variability for several classes of univalent functions has become an interesting area of research in geometric function theory. Several authors have studied the region of variability for several classes of functions. In $2005$, Yanagihara \cite{yanagihara} studied region of variability for functions of bounded derivatives. In $2006$, Yanagihara \cite{yanagihara2} proved the region of variability for convex functions. In $2007$, Ponnusamy and Vasudevarao \cite{Ponnusamy-Allu-2007-JMMA} studied the region of variability of two subclases of univalent functions. In $2009$, Ponnusamy {\it et al.} \cite{S.Ponnusamy variability regions-2009} have studied the region of variability for the class $\mathcal{B}(\alpha, \beta,M)$ of analytic and univalent funtions in the unit disk $\mathbb{D}$ with $f(0)=0$, $f'(0)=\alpha$, and $f''(0)=M\beta$ satisfying $|zf''(z)|\leq M$ for $\alpha, \beta \in \mathbb{C}$ and $M \in \mathbb{R}$ with $0<M \leq |\alpha|$ and $|\beta|\leq 1$. Furthermore, Ponnusamy {\it et al.} \cite{S.Ponnusamy variability regions-2009} have also studied the region of variability for the class $\mathcal{P}(\alpha,M)$ of analytic and univalent funtions in $\mathbb{D}$ with $f(0)=0$, $f'(0)=\alpha$, and $f''(0)=M\beta$ satisfying $\real zf''(z)> -M$ in $\mathbb{D}$, where $\alpha \in \mathbb{C} \setminus \{0\}$ and $0<M \leq 1/\log\, 4$. In $2009$,  Ponnusamy {\it et al.} \cite{Ponnu-Allu-2009-Colloq} studied the region of variability for spiral-like functions with respect to a boundary point. In $2010$, Ponnusamy and Vasudevarao \cite{Ponnusamy-Allu-2010-Polo} extensively studied the region of variability for functions with positive real part. In $2011$, Ponnusamy {\it et al.} \cite{Ponnu-Allu-2011-CAOT} investigated the region of variability for exponentially convex univalent functions. For other interesting results on region of variability, we refer to \cite{Ponnusamy-Allu-2007-JMMA,Ponnusamy-Allu-2008-CVEE,S.Ponnusamy variability regions-2013}.\\

If $f$ is locally univalent function, then pre-Schwarzian norm is defined by 
\begin{equation}
	\|f\|= \sup_{z\in \mathbb{D}}\, (1-|z|^2)\bigg|\dfrac{f''(z)}{f'(z)}\bigg|.
\end{equation}
The class of univalent functions is preserved under a number of elementary transformations. The preservation of $\mathcal{S}$ under the disk automorphism leads to the study of the behaviour of the pre-Schwarzian norm of $f$. The pre-Schwarzian norm has the close connection with functions in $\mathcal{S}$. Hence the study of region of variability problems together with pre-Schwarzian norm is very classic and worthy. It is well known that $||f|| \leq 6$ if f is univalent in $\mathbb{D}$, and conversely if $||f|| \leq 1$ then $f$ is univalent in $\mathbb{D}$, and these bounds are sharp (see \cite{becker-1984}).

\vspace{2mm}
The following lemma is useful to prove our main result.

\begin{lem}\label{lem-2.11} \cite{Ponnusamy-Allu-2007-JMMA}
	For $ \theta\in\mathbb{R} $ and $ 0<\lambda<1 $, the function\begin{align*}
		G(z)=\int_{0}^{z}\frac{e^{i\theta}z}{\left(1+\left(e^{i\theta}-1\right)\lambda z-e^{i\theta}z\right)^2} \,dz,\; z\in\mathbb{D}
	\end{align*}
	has a double zero at the origin and no zeros elsewhere in $ \mathbb{D} $. Furthermore, these exists a starlike univalent function $ G_0 $ in $ \mathbb{D} $ such that $ G=2^{-1}e^{i\theta}G^2_0 $ and $ G_0(0)=G^{\prime}_0(0)-1 $.
\end{lem}
For a positive integer $p$, let 
\begin{equation*}
	\left(\mathcal{S}^{*}\right)^p= \{f= f_{0}^p: f_{0} \in \mathcal{S}^*\}.
\end{equation*}
\begin{lem} \label{lem-a-p-1}\cite{yanagihara}
Let $f$ be an analytic function in $\mathbb{D}$ with $f(z)=z^p + \cdots .$ If
\begin{equation*}
\real \left(1+ z\frac{f''(z)}{f'(z)}\right)>0, \,\,\, z \in \mathbb{D},
\end{equation*}
then $f \in \left(\mathcal{S}^{*}\right)^p$.
\end{lem}
\section{Preliminaries}
For $0<\alpha<1$ and $0<\beta <1$, let 
\begin{equation} \label{e-p-2.1}
	f_{\alpha}(z)=\displaystyle\frac{z}{1-z^{\alpha}}, \,\, z \in \mathbb{D}
\end{equation}
and
\begin{align*}
	\mathcal{F}=\bigg\{f_{\alpha} \in \mathcal{A}: {\rm Re}\left(1+\frac{zf^{\prime\prime}_{\alpha}(z)}{f^{\prime}_{\alpha}(z)}\right)>\beta\bigg\}.
	\end{align*}
Let ${\mathcal F}(\lambda)  =\{f_{\alpha}\in \mathcal{F}:\,f_{\alpha}''(0)=2\lambda (1-\beta) \,\,\, \mbox{for} \,\,\, 0\leq \lambda \leq 1\}$.
For $f_{\alpha} \in \mathcal{F}$, let 
\begin{equation}\label{e-2.1-b}
	P_f(z)=1+ \displaystyle \frac{zf^{\prime\prime}_{\alpha}(z)}{f^{\prime}_{\alpha}(z)}-\beta.
\end{equation}
Clearly, $ P_f(0)=1-\beta $. Let $\mathcal{B}_{0}$ be the class of analytic functions $\omega$ in the unit disk $\mathbb{D}$ such that $|\omega(z)|\leq 1$ in $\mathbb{D}$ and $\omega(0)=0$. Then for each $f \in \mathcal{F}(\lambda)$, there exists an $\omega_f \in \mathcal{B}_{0}$ of the form 
\begin{align} \label{e-2.1-a}
	\omega_f(z):=\frac{P_f(z)-(1-\beta)}{P_f(z)+(1-\beta)}.
\end{align}
A simple computation shows that 
\begin{align*}
	\omega^{\prime}_f(0)=\lim_{z\rightarrow 0}\frac{\omega_f(z)-\omega_f(0)}{z-0}=\lim_{z\rightarrow 0}\frac{\omega_f(z)}{z}.
\end{align*}
Let $ \omega_f^{\prime}(0)=\lambda $, where $ |\lambda|\leq 1 $. From \eqref{e-2.1-a}, we have 
\begin{equation*}
\omega_f^{\prime}(z)= \dfrac{2(1-\beta)P'_{f}(z)}{\left(P_{f}(z) + (1-\beta)\right)^2}.	
\end{equation*}
Since $P_{f}(0)=1-\beta$, we have 
\begin{equation} \label{e-p-2}
\omega^{\prime}_f(0)= \frac{P'_{f}(0)}{2(1-\beta)}.	
\end{equation}
Therefore
\begin{equation} \label{e-p-2-a}
	\lambda= \frac{P'_{f}(0)}{2(1-\beta)}.
\end{equation}
From \eqref{e-2.1-b}, we have 
\begin{equation} \label{e-p-3}
P_{f}(z) f'_{\alpha}(z)= (1-\beta) f'_{\alpha}(z) + z f''_{\alpha}(z).	
\end{equation}
Differentiating \eqref{e-p-3} gives 
\begin{equation} \label{e-p-4}
P'_{f}(z) f'_{\alpha}(z) + \left(P_{f}(z) - (1-\beta) - 1\right) f''_{\alpha}(z) + zf''_{\alpha}(z)=0.
\end{equation}
A simple computation of \eqref{e-p-4} at $z=0$ shows that
\begin{equation} \label{e-p-5}
P'_{f}(0)=f''_{\alpha}(0).
\end{equation}
From \eqref{e-p-2-a} and \eqref{e-p-5}, we obtain
\begin{equation}
	f''_{\alpha}(0)=2\lambda(1-\beta).
\end{equation}
Define $ g : \mathbb{D}\rightarrow\mathbb{D} $ by $$ g(z)=\displaystyle\frac{\frac{\omega_f(z)}{z}-\lambda}{1-\frac{\overline{\lambda}\omega_f(z)}{z}},
$$
equivalently
$$
g(z)=\dfrac{\omega_f(z) - \lambda z}{z-\overline{\lambda}\omega_f(z)}.
$$
Then, clearly $ g(0)=0 $ and $ g $ is an analytic function in $\mathbb{D}$. Therefore $g \in \mathcal{B}_{0}$.
\vspace{1.2mm}

First we prove that the class $\mathcal{F}(\lambda)$ is rotationally invariant.
The class $\mathcal{F}(\lambda)$ is rotationally invariant if, and only if, $f_{\theta} \in \mathcal{F}(\lambda)$, where $f_{\theta}:=e^{-i\theta}f_{\alpha}(e^{i\theta}z)$, where $f_{\alpha}$ is defined by \eqref{e-p-2.1}.
A simple computation shows that
\begin{align*}
	f^{\prime}_{\theta}(z)=f_{\alpha}^{\prime}\left(e^{i\theta}z\right)\;\; \mbox{and} \,\, f^{\prime\prime}_{\theta}(z)=e^{i\theta}f_{\alpha}^{\prime\prime}\left(e^{i\theta}z\right).
\end{align*}
Then, it is easy to see that
\begin{align*}
	{\rm Re}\left(1+\frac{zf^{\prime\prime}_{\theta}(z)}{f^{\prime}_{\theta}(z)}\right)={\rm Re}\left(1+\frac{ze^{i\theta}f^{\prime\prime}\left(e^{i\theta}z\right)}{f^{\prime}\left(e^{i\theta}z\right)}\right)={\rm Re}\left(1+\frac{\xi f''(\xi)}{f'(\xi)}\right)>\beta,
\end{align*}
where $ \xi=ze^{i\theta} $. Since $f_{\alpha} \in \mathcal{F}(\lambda)$, it is easy to see that $f_{\theta}(0)=0$, $f'_{\theta}(0)=1$, and $f''_{\theta}(0)=2\lambda (1-\beta) e^{i \theta}$. Therefore, $f_{\theta} \in \mathcal{F}(\lambda)$. Hence the class $ \mathcal{F}(\lambda) $ is rotational invariant. Since the class $\mathcal{F}(\lambda)$ is rotationally invariant, without loss of generality we may assume that $0\leq \lambda \leq 1$ instead of $\lambda \in \mathbb{D}$.
\\

We wish to construct an extremal function $F_{a,\lambda}(z)$ where $a \in \overline{\mathbb{D}}=\{z\in \mathbb{C}:|z|\leq 1\}$ and it belongs to the class $\mathcal{F}(\lambda)$.
From the second condition of the Schwarz lemma, we have $ g(z)=az $, where $ |a|=1 $. Therefore, we see that
\begin{align*}
g(z)=\dfrac{\omega_f(z) - \lambda z}{z-\overline{\lambda}\omega_f(z)}=	\frac{P_f(z)-(1-\beta)\left(\dfrac{1+\lambda z}{1-\lambda z}\right)}{P_f(z)+(1-\beta)\left(\dfrac{z+\lambda}{z-\lambda}\right)}=az
\end{align*}
which implies that
\begin{align}\label{e-2.4}
	\dfrac{\displaystyle\frac{f_{\alpha}^{\prime\prime}(z)}{f_{\alpha}^{\prime}(z)}-\frac{2(1-\beta)\lambda}{1-\lambda z}}{\displaystyle\frac{f_{\alpha}^{\prime\prime}(z)}{f_{\alpha}^{\prime}(z)}+\frac{2(1-\beta)}{z-\lambda}}=az\left(\frac{z-\lambda}{1-\lambda z}\right):=az \, \delta(z,\lambda), \,\,\mbox{where} \,\,\delta(z,\lambda)= \frac{z-\lambda}{1-\lambda z}.
\end{align}
Further simplification of \eqref{e-2.4} gives
\begin{align}\label{e-2.4-a}
	\dfrac{f_{\alpha}^{\prime\prime}(z)}{f_{\alpha}^{\prime}(z)}=2(1-\beta)\dfrac{\dfrac{\lambda}{1-\lambda}+\dfrac{az\delta(z, \lambda)}{z-\lambda}}{1-az\delta(z, \lambda)}=\dfrac{2(1-\beta)(az+\lambda)}{1-\lambda z-az^2 + \lambda a z}.
\end{align}
From \eqref{e-2.4-a}, we obtain 
\begin{align}\label{e-2.5}
	\frac{f_{\alpha}^{\prime\prime}(z)}{f_{\alpha}^{\prime}(z)}=2(1-\beta)\frac{\delta(az, \lambda)}{1-\delta(az,\lambda)z}.
\end{align}
Integrating both sides of \eqref{e-2.5}, we obtain
\begin{align} \label{e-2.7-b}
	F_{a,\lambda}(z):=f_{\alpha}(z)=\int_{0}^{z}\exp\left(\int_{0}^{\xi}\frac{2(1-\beta)\delta(a\xi,\lambda)}{1-\delta(a\xi,\lambda)}\,\,d\xi\right)d\xi.
\end{align}
Clearly, $F_{a,\lambda}(0)=0, F_{a,\lambda}'(0)=1$ and $F_{a,\lambda}''(0)=2\lambda (1-\beta)$.
Hence $	F_{a,\lambda} \in \mathcal{F}(\lambda)$.
\vspace{2mm}

If $f\in \mathcal{F}(\lambda)$, we define 
\begin{equation} \label{e-2.27-b}
	h(z)=\frac{P_{f}(z)}{1-\beta}.
\end{equation}
where 
$P_{f}(z)=1+ zf_{\alpha}''/{f_{\alpha}'}-\beta.$ 
 So that ${\rm Re\,} h(z) > 0$ and $h(0)=1$.
\vspace{2mm}
Therefore by Herglotz representation for functions with positive real part there exists unit probability measure $\mu$ on $(-\pi,\pi]$, {\it{i.e.,}}
$$ \int_{-\pi}^{\pi}\,d\mu (t)=1,$$
such that 
\begin{equation} \label{e-h-2.25}
	h(z)=\int_{-\pi}^{\pi}\frac{1+ze^{-it}}{1-ze^{-it}}\,\,d\mu(t).
\end{equation}
In view of \eqref{e-2.27-b}, we have 
$$P_{f}(z)=(1-\beta)\int_{-\pi}^{\pi}\frac{1+ze^{-it}}{1-ze^{-it}}\,\,d\mu(t). $$
From \eqref{e-0.1}, we obtain
\begin{equation}
	1+ \frac{zf_{\alpha}''}{f_{\alpha}'}-\beta= (1-\beta)\int_{-\pi}^{\pi}\frac{1+ze^{-it}}{1-ze^{-it}}\,\,d\mu(t),
\end{equation}
which implies that 
$$ \frac{zf_{\alpha}''}{f_{\alpha}'} = (1-\beta)\int_{-\pi}^{\pi}\bigg(\frac{1+ze^{-it}}{1-ze^{-it}} -1\bigg)\,\,d\mu (t).$$
Therefore,
$$\frac{f_{\alpha}''}{f_{\alpha}'} =2(1-\beta) \int_{-\pi}^{\pi}\frac{e^{-it}}{1-ze^{-it}}\,\,d\mu(t).
$$
Integrating on both sides with respect to $z$, we obtain
$$\log f_{\alpha}'(z) =2(1-\beta) \int_{-\pi}^{\pi}-\log (1-ze^{-it})\,\,d\mu (t),$$
and hence
\begin{equation}\label{eq-0}
	f_{\alpha}(z)=\int_{0}^{t}\exp \bigg(2(1-\beta) \int_{-\pi}^{\pi}-\log (1-ze^{-it})\, d\mu (t)\bigg).
\end{equation}

The main aim of this paper is to obtain the region of variability of $\log f'(z_0)$ where $f$ runs over the class $\mathcal{F}(\lambda)$. Further, we demonstrate the region of variability of $\log f'(z_0)$ for the class $\mathcal{F}(\lambda)$ by using Mathematica.

\section{Main Results}

Consider the following set 
$$V(z_0, \lambda, \beta)=\{\log f^{\prime}(z_0): f\in\mathcal{F}(\lambda)\} .$$
Now we mention some basic properties of the set $V(z_0, \lambda, \beta)$.
\begin{prop}\label{prop-5.1}
	We have
	\begin{itemize}
		\item[(1)]  $V(z_0, \lambda, \beta)$ is compact.
		\item[(2)] $V(z_0, \lambda, \beta)$ is convex.
		\item[(3)] For $|\lambda|=1$ or $z_{0}=0$, we have
		\begin{equation}\label{vvv}
		V(z_0, \lambda, \beta) = \{-2(1-\beta) \log (1-\lambda z)\}	
		\end{equation}
		\item[(4)] For $|\lambda| <1$ and $z_{0} \in \mathbb{D} \setminus \{0\}$, $V(z_0, \lambda, \beta)$ has $-2(1-\beta) \log (1-\lambda z)$ as an interior point.
	\end{itemize}
\end{prop}

	\begin{prop} \label{prop-1}
		For $ f\in \mathcal{F}(\lambda) $ where $ 0\leq \lambda<1 $, we have
		\begin{align*}
			\left|\frac{f^{\prime\prime}(z)}{f^{\prime}(z)}-c(z,\lambda,\beta)\right|\leq r(z,\lambda,\beta),\; z\in\mathbb{D}.
		\end{align*}
where
$$
c(z,\lambda, \beta)=\frac{2(1-\beta)\left(\lambda(1-|z|^2)+(|z|^2-\lambda^2)\right)\bar{z}}{(1-|z|^2)(1-2\lambda({\rm Re}z)+|z|^2)}
$$
and
$$ r(z,\lambda, \beta)=
	\frac{2(1-\beta)(1-\lambda^2)|z|}{(1-|z|^2)(1-2\lambda({\rm Re}z)+|z|^2)}
$$		
For each $z \in \mathbb{D}\setminus\{0\}$,equality hold if, and only if, $f=F_{e^{i\theta},\lambda}$ for some $\theta \in \mathbb{R}$.
	\end{prop}

\begin{cor} \label{cor-3.6}
	Let $ \gamma :=z(t)  $, where $ 0\leq t\leq 1 $ be a $C^{1}$ curve in $ \mathbb{D} $ with $ z(0)=0 $ and $ z(1)=z_0 $. Then, we have
	\begin{align*}
		V(z_0, \lambda, \beta)\subset \bigg\{w\in\mathbb{C} : |\omega-C(\lambda,\gamma,\beta)|\bigg\}\leq R(\lambda,\gamma,\beta),
	\end{align*}
	where $ C(\lambda,\gamma,\beta)=\displaystyle\int_{0}^{1}c(z(t),\lambda,\beta)z^{\prime}(t)dt $ and 
	$ R(\lambda,\gamma, \beta)=\displaystyle\int_{0}^{1}r(z(t),\lambda,\beta)|z^{\prime}(t)|dt $.
\end{cor}

Now we state the following result which is vital to prove our main result.
\begin{prop}(Uniqueness of the curve) \label{prop-2}
	Let $ z_0\in\mathbb{D}\setminus\{0\}$. Then, for $ \theta\in (-\pi, \pi] $, we have \begin{align*}
		\log F^{\prime}_{e^{i\theta},\lambda}(z_0)\in \partial V(z_0,\lambda, \beta).
	\end{align*}
Furthermore, if $ \log f^{\prime}(z_0)=\log F^{\prime}_{e^{i\theta},\lambda}(z_0)  $ for some $ f\in\mathcal{F}(\lambda) $ and $\theta\in (-\pi, \pi]$, then $ f=F_{e^{i\theta},\lambda} $. 
\end{prop}
Finally, we state our main result to obtain the region of variability of $\log f'(z_0)$ where $f$ runs over the class $\mathcal{F}(\lambda)$.
\bthm\label{pvdev3-th1}
For $0\leq\lambda <1$ and $z_0\in\mathbb{D}\setminus \{0\}$,
the boundary $\partial{V(z_0,\lambda,\beta)}$ is the Jordan curve given
by
\be\label{pvdev3-eq7}
(-\pi,\pi]\ni \theta \mapsto \log F'_{e^{i\theta},\lambda}(z_0)
 = \int_0^{z_0} \frac{\delta(e^{i\theta}z,\lambda)}
{z\delta(e^{i\theta}z,\lambda)-1}\,dz, \quad z\in\mathbb{D}.
\ee
If  $\log f'(z_0)=\log F'_{e^{i\theta},\lambda}(z_0)$ for some
$f\in \mathcal{F}(\lambda)$ and $\theta\in (-\pi,\pi]$, then $f(z)=
F_{e^{i\theta},\lambda}(z)$. Here $F_{e^{i\theta},\lambda}(z)$ is
given by 
\begin{equation}\label{ppp}
	F_{e^{i\theta},\lambda}(z)= f_{\alpha}(z)=\int_{0}^{z}\exp\left(\int_{0}^{\xi}\frac{2(1-\beta)\delta(e^{i\theta}\xi,\lambda)}{1-\delta(e^{i\theta}\xi,\lambda)}d\xi\right)d\xi.
\end{equation}
\ethm


\section{Proof of the main results}

\begin{pf}[{\bf Proof of Proposition \ref{prop-5.1}}]
	(1) To show that $V(z_0, \lambda, \beta)$ we need to show that the set $\mathcal{F}(\lambda)$ is compact. For this, we need to show that $\mathcal{F}(\lambda)$ is closed and locally bounded. Let $f_{n}$ be a sequence of functions in $F(\lambda)$ such that $f_{n} \rightarrow f$ uniformly on every compact subsets of the unit disk. Then to show that $F(\lambda)$ is compact, we show that $f \in F(\lambda)$. Since $f_{n} \rightarrow f$ uniformly on every compact subset of $\mathbb{D}$, we have $f'_{n} \rightarrow f'$ uniformly on every compact subsets of $\mathbb{D}$. Therefore, if $f_{n} \rightarrow f$ uniformly, then $f_{n}(0) \rightarrow f(0)$ pointwise which gives that $f(0)=0$. Since, $f'_{n} \rightarrow f'$ uniformly on every compact subsets of $\mathbb{D}$, we get
\begin{equation} \label{e-a-p-1}
	f'_{n}(0) \rightarrow f'(0), \,\, f''_{n}(0) \rightarrow f''(0)
\end{equation}
pointwise. Since $f_{n} \in \mathcal{F}(\lambda)$, $f'_{n}(0)=1$ and $f''_{n}(0)=2\lambda(1-\beta)$. In view of \eqref{e-a-p-1}, we have  $f'(0)=1$ and $f''(0)=2\lambda(1-\beta)$. In view of the above discussion, it is easy to see that 
\begin{equation*}
	1+ \dfrac{z f''_{n}(z)}{f'_{n}(z)} \longrightarrow 1+ \dfrac{z f''(z)}{f'(z)}
\end{equation*}
uniformly on every compact subsets of $\mathbb{D}$. Since $f_{n} \in \mathcal{F}(\lambda)$, we have 
\begin{equation*}
	0<\real \left(1+ \dfrac{z f''_{n}(z)}{f'_{n}(z)}\right) \longrightarrow \real \left(1+ \dfrac{z f''(z)}{f'(z)}\right) \geq 0
\end{equation*}
in $\mathbb{D}$. Therefore, by the maximum modulus principle, it follows that 
\begin{equation*}
\real \left(1+ \dfrac{z f''(z)}{f'(z)}\right) >0	
\end{equation*}
in $\mathbb{D}$. 	
Locally boundedness follows from the definition itself, thus $\mathcal{F}(\lambda)$ is compact and in view of the continuous function $\psi:\mathcal{F}(\lambda)\to V(z_0, \lambda, \beta)$ defined by 
	$$\psi(f)=log(f'(z_0))$$
	we conclude that the class $V(z_0, \lambda, \beta)$ is also compact.\\
	(2) Let $f_{1}, f_{2} \in \mathcal{F}(\lambda)$ and $0 \leq t \leq 1$. Then it is easy to see that the function 
	\begin{equation}
		f_{t}(z)= \int_{0}^{z} \exp \left\{(1-t) \log f'_{1}(\xi) + t \log f'_{2}(\xi)\right\} \, d\xi
	\end{equation}
	belongs to $\mathcal{F}(\lambda)$ and hence $V(z_{0}, \lambda, \beta)$ is convex.\\
	(3) If $z_{0}=0$, then \eqref{vvv} is trivially true. If $|\lambda|=|\omega'_{f}(0)|=1$, then by an application of the Schwarz lemma, we have $\omega_f(z)=\lambda z$. Therefore using \eqref{e-2.1-a}, we obtain 
	\begin{equation} \label{e-a-p-2}
		P_{f}(z)= (1-\beta) \left(\frac{1+\lambda z}{1-\lambda z}\right).
	\end{equation} 
	Further simplification of \eqref{e-a-p-2} by using \eqref{e-2.1-b} shows that 
	\begin{equation*}
		\log f'(z)= 2(1-\beta) \frac{\lambda}{1-\lambda z}
	\end{equation*}
	which implies that 
	\begin{equation*}
		f(z)=\int_{0}^{z} \exp\left(2(1-\beta) \frac{\lambda}{1-\lambda \xi}\right) \, d\xi.
	\end{equation*}
	Conversely, 
	\begin{equation*}
		V(z_{0},\lambda,\beta)=\{-2(1-\beta) \log(1-\lambda z)\}.
	\end{equation*}

(4) For $|\lambda|<1$ and $a \in \overline{\mathbb{D}}$, we define 
\begin{equation*}
	\delta(z, \lambda)=\frac{z+\lambda}{1+\overline{\lambda}z}
\end{equation*}
and 
\begin{equation}\label{ppp}
	F_{a,\lambda}(z)= f_{\alpha}(z)=\int_{0}^{z}\exp\left(\int_{0}^{\xi}\frac{2(1-\beta)\delta(a\xi,\lambda)}{1-\delta(a\xi,\lambda)}d\xi\right)d\xi.
\end{equation}
An observation shows that 
\begin{equation}\label{main}
	\omega_{F_{a,\lambda}}(z)=z \delta (az,\lambda). 
\end{equation}
We want to prove that the mapping $\mathbb{D} \ni a \mapsto \log F'_{a, \lambda}(z_{0})$ is a non-constant analytic function of "$a$" for each fixed $z_{0} \in \mathbb{D} \setminus \{0\}$ and $\lambda \in \mathbb{D}$. To show this, we let 
\begin{align}
	h(z)&= \left .\frac{1}{(1-\lambda ^2)}
	\frac{\partial}{\partial a}\left\{\frac{}{}\log F'_{a,\lambda}(z)\right\}\right
	|_{a=0} \\ \nonumber
	&= 2(1-\beta) \int_{0}^{z} \frac{\xi}{(1-\lambda \xi )^2} \,\, d\xi \\ \nonumber
	&=z^2 + \cdots .
\end{align}

A simple computation shows that 
\begin{equation*}
	\real \left\{\frac{zh''(z)}{h'(z)}\right\}=2(1-\beta) \real \left\{\frac{1+\lambda z}{1-\lambda z}\right\}>0 \,\, \mbox{for} \,\, z\in \mathbb{D}.
\end{equation*}
By Lemma \ref{lem-a-p-1}, there exists a function $h_{0} \in \mathcal{S}^{*}$ with $h=h_{0}^2$. The univalence of $h_{0}$ together with the condition $h_{0}(0)=0$ implies that $h(z_{0})\neq 0$ for $z_{0} \in \mathbb{D} \setminus \{0\}$. Consequently, the mapping $\mathbb{D} \ni a \mapsto \log F'_{a, \lambda}(z_{0})$ is a non-constant analytic function of "$a$", and hence it is an open mapping. Therefore $V(z_{0},\lambda,\beta)$ contains the open set $\{\log F'_{a, \lambda}(z_{0}): |a|<1 \}$. In particular, 
\begin{equation*}
	\log F'_{0, \lambda}(z_{0})= -2(1-\beta) \log (1-\lambda z)
\end{equation*}
is an interior point of $\{\log F'_{a, \lambda}(z_{0}): |a|<1 \} \subset V(z_{0},\lambda,\beta)$.
This completes the proof.
\end{pf}
\begin{pf}[{\bf Proof of Proposition \ref{prop-1}}]
	For $f \in \mathcal{F}(\lambda)$, define $g(z)$ by  
	$$
	g(z)=\frac{	\omega_f(z)-\lambda z}{z-\overline{\lambda}\omega_f(z)}.
	$$
where $w_f$ is defined by \eqref{e-2.1-a}.
\\

From \eqref{e-2.1-a}, it is easy to see that 
	\begin{align}\label{e-2.1}
		\omega_f(z)-\lambda z=\frac{P_f(z)\left(1-\lambda z\right)-(1-\beta)(1+2\lambda)}{P_f(z)+(1-\beta)}
	\end{align}
	and 
	\begin{align}\label{e-2.2}
		z-\overline{\lambda}\omega_f(z)=\frac{(z-\lambda)P_f(z)+(1-\beta)(z+\lambda)}{P_f(z)+(1-\beta)}.
	\end{align}
	Therefore, in view \eqref{e-2.1} and \eqref{e-2.2}, we obtain
	\begin{align*}
		g(z)=\frac{	\omega_f(z)-\lambda z}{z-\overline{\lambda}\omega_f(z)}=\frac{P_f(z)(1-\lambda z)-(1-\beta)(1+\lambda z)}{(z-\lambda)P_f(z)+(1-\beta)(z+\lambda)}.
	\end{align*}
	Again, $ |g(z)|\leq |z| $ shows that
	\begin{align*}
		\ds	\left|\frac{P_f(z)(1-\lambda z)-(1-\beta)(1+\lambda z)}{(z-\lambda)P_f(z)+(1-\beta)(z+\lambda)}\right|\leq |z|
	\end{align*}
	which can be written as 
	\begin{align}\label{e-2.7-a}
		\left|\frac{P_f(z)-(1-\beta)\left(\dfrac{1+\lambda z}{1-\lambda z}\right)}{P_f(z)+(1-\beta)\left(\dfrac{z+\lambda}{z-\lambda}\right)}\right|\leq |z|\left|\dfrac{z-\lambda}{1-\lambda z}\right|.
	\end{align}
	Further simplification of \eqref{e-2.7-a} by using \eqref{e-2.1-b}, we obtain
	\begin{align*}
		\left|\frac{(1-\beta)+\dfrac{zf^{\prime\prime}_{\alpha}(z)}{f_{\alpha}^{\prime}(z)}-(1-\beta)\left(\dfrac{1+\lambda z}{1-\lambda z}\right)}{(1-\beta)+\dfrac{zf^{\prime\prime}_{\alpha}(z)}{f_{\alpha}^{\prime}(z)}+(1-\beta)\left(\dfrac{z+\lambda}{z-\lambda}\right)}\right|\leq |z|\left|\frac{z-\lambda}{1-\lambda z}\right|
	\end{align*}
	which implies that
	\begin{align} \label{e-2.9-a}
		\left|\dfrac{\dfrac{zf^{\prime\prime}_{\alpha}(z)}{f_{\alpha}^{\prime}(z)}-(1-\beta)\left(\dfrac{2\lambda z}{1-\lambda z}\right)}{\dfrac{zf^{\prime\prime}_{\alpha}(z)}{f_{\alpha}^{\prime}(z)}+(1-\beta)\left(\frac{2z}{z-\lambda}\right)}\right|\leq |z|\left|\dfrac{z-\lambda}{1-\lambda z}\right|.	
	\end{align}
	A simple computation shows that \eqref{e-2.9-a} is equivalent to 
	\begin{align}\label{e-2.3}
		\left|\displaystyle\dfrac{\dfrac{zf^{\prime\prime}_{\alpha}(z)}{f_{\alpha}^{\prime}(z)}-A(z,\lambda,\beta)}{\displaystyle\dfrac{zf^{\prime\prime}_{\alpha}(z)}{f_{\alpha}^{\prime}(z)}+B(z,\lambda,\beta)}\right|\leq |z|\left|\tau\left(z,\lambda\right)\right|,
	\end{align}
	where 
	$$ 
	A(z,\lambda, \beta)= \frac{2(1-\beta)\lambda}{1-\lambda z}, B(z, \lambda, \beta)= \frac{2(1-\beta)z}{z-\lambda},\,\, \mbox{and}\,\, \tau(z,\lambda)=\frac{z-\lambda}{1-\lambda z}.
	$$ 
	Further simplification of \eqref{e-2.3} implies that
	\begin{align*}
		&\left|\frac{f^{\prime\prime}_{\alpha}(z)}{f_{\alpha}^{\prime}(z)}-\frac{A(z,\lambda,\beta)+|z|^2 |\tau(z,\lambda)|^2B(z,\lambda, \beta)}{1-|z|^2|\tau(z,\lambda)|^2}\right|  \\[3mm]
		& \leq \dfrac{|z| |\tau(z,\lambda)| |A(z,\lambda,\beta) + B(z, \lambda)|}{1-|z|^2|\tau(z,\lambda)|^2},
	\end{align*}
	where \begin{align*}
		A(z,\lambda,\beta)+|z|^2 |\tau(z,\lambda)|^2B(z,\lambda, \beta)&=2(1-\beta)\left(\frac{\lambda}{1-\lambda z}+|z|^2\frac{(z-\lambda)(\bar{z}-\lambda)}{|1-\lambda z|^2}\frac{1}{z-\lambda}\right)\\&=2(1-\beta)\left(\frac{\lambda(1-|z|^2)+\left(|z|^2-\lambda^2\right)\bar{z}}{|1-\lambda z|^2}\right).
	\end{align*}
	A simple computation shows that
	\begin{align*}
		A(z,\lambda,\beta)+B(z,\lambda,\beta)=2(1-\beta)\frac{1-\lambda^2}{(1-\lambda z)(z-\lambda)}.
	\end{align*}
	Therefore, it is easy to see that
	\begin{align*}
		|z||\tau(z,\lambda)| |A(z,\lambda,\beta)+B(z,\lambda,\beta)|=\frac{2(1-\beta)(1-\lambda^2)|z|}{|1-\lambda z|^2}
	\end{align*}
	which shows that
	\begin{align*}
		1-|z|^2|\tau(z,\lambda)|^2=\frac{(1-|z|^2)(1-2\lambda \real z)+|z|^2}{|1-\lambda z|^2}.
	\end{align*}
	Thus, we have
	\begin{align*}
		\frac{A(z,\lambda,\beta)+|z|^2|\tau(z,\lambda)B(z,\lambda,\beta)|}{1-|z|^2|\tau(z,\lambda)|}&=\frac{2(1-\beta)\left(\lambda(1-|z|^2)+(|z|^2-\lambda^2)\right)\bar{z}}{(1-|z|^2)(1-2\lambda({\rm Re}z)+|z|^2)}\\&:=c(z,\lambda, \beta).
	\end{align*}
Further,
\begin{align*}
	\cfrac{|z||\tau(z,\lambda)| |A(z,\lambda,\beta)+B(z,\lambda,\beta)|}{1-|z|^2|\tau(z,\lambda)|} &=
	\frac{2(1-\beta)(1-\lambda^2)|z|}{(1-|z|^2)(1-2\lambda({\rm Re}z)+|z|^2)}\\&:=r(z,\lambda, \beta).
	\end{align*}
	This completes the proof.
\end{pf}

We note that if $f \in \mathcal{F}(0)$, then pre-Schwarzian norm 
$\|f\|\leq 2(1-\beta).$ Indeed, 
\begin{equation}\label{e-0.3}
	\bigg|\frac{f''(z)}{f'(z)}-c(z,0,\beta)\bigg|\leq r(z,0,\beta),
\end{equation}
where
$$ c(z,0,\beta)= \dfrac{2(1-\beta)|z|^2\bar{z}}{(1-|z|^2)(1+|z|^2)}\,\,\,\, \mbox{and}\,\,\,\, r(z,0, \beta)=\dfrac{2(1-\beta)|z|}{(1-|z|^2)(1+|z|^2)}.$$
From \eqref{e-0.3}, we have
$$\bigg|\dfrac{f''(z)}{f'(z)}- \dfrac{2(1-\beta)|z|^2\bar{z}}{(1-|z|^2)(1+|z|^2)}\bigg| \leq \dfrac{2(1-\beta)|z|}{(1-|z|^2)(1+|z|^2)}$$
which implies
$$\bigg|\dfrac{f''(z)}{f'(z)}\bigg|- \bigg|\dfrac{2(1-\beta)|z|^2\bar{z}}{1-|z|^4}\bigg|\leq \dfrac{2(1-\beta)|z|}{1-|z|^4}.
$$
Therefore,
$$ \bigg|\dfrac{f''(z)}{f'(z)}\bigg| \leq \dfrac{2(1-\beta)|z|}{1-|z|^4} + \dfrac{2(1-\beta)|z|^3}{1-|z|^4} = \dfrac{2(1-\beta)|z|}{1-|z|^2}
$$
which implies that 
$$(1-|z|^2)\bigg|\dfrac{f''(z)}{f'(z)}\bigg| \leq 2(1-\beta)|z|, $$
and hence
$$\sup_{z\in \mathbb{D}} \,(1-|z|^2)\bigg|\dfrac{f''(z)}{f'(z)}\bigg| \leq 2(1-\beta) \,\, \text{for} \,\, z\in \mathbb{D}.
$$
Therefore, if $f\in \mathcal{F}(0)$, then 
$$\|f\|\leq 2(1-\beta). $$ 

\begin{proof}[{\bf Proof of Corollary \ref{cor-3.6}}]
	A simple computation using the proposition \ref{prop-1} shows that
	\begin{align*}
		\left|\log f^{\prime}_{\alpha}(z)-C(\lambda,\gamma,\beta)\right|&=\left|\log f^{\prime}_{\alpha}(z)-\int_{0}^{1}c(z(t),\lambda,\beta)z^{\prime}(t)\, dt \right|\\  & = \left|\bigintsss_{0}^{1}\frac{f^{\prime\prime}_{\alpha}(z(t))}{f^{\prime}_{\alpha}(z(t))}z^{\prime}(t)dt-\bigintsss_{0}^{1}c(z(t),\lambda,\beta)z^{\prime}(t)\,dt\right|\\&  \leq \bigintsss_{0}^{1}\left|\left(\frac{f^{\prime\prime}_{\alpha}(z(t))}{f^{\prime}_{\alpha}(z(t))}-c(z(t),\lambda,\beta)\right)\right||z^{\prime}(t)|\, dt\\&
		\leq \int_{0}^{1}r(z(t),\lambda,\beta)|z^{\prime}(t)|\,dt:=R(\lambda,\gamma,\beta).
	\end{align*}
	This completes the proof.
\end{proof}	

\begin{pf}[{\bf Proof of Proposition \ref{prop-2}}]
	From \eqref{e-2.7-b}, we have 
	$$ F_{a,\lambda}(z)=\displaystyle\int_{0}^{z}\exp\left(\int_{0}^{\xi_{2}}\frac{2(1-\beta)\delta(a\xi_{1},\lambda)}{1-\delta(a\xi_{1},\lambda)}\,\,d\xi_{1}\right)d\xi_{2} .
	$$ 
	Then, it is easy to see that
	\begin{align*}
		F^{\prime}_{e^{i \theta},\lambda}=\exp\left(\int_{0}^{z}\frac{2(1-\beta)\delta(a\xi_{1},\lambda)}{1-\delta(a\xi_{1},\lambda)}\,\, d\xi_{1}\right)
	\end{align*}
	and hence,
	\begin{align*}
		\frac{F^{\prime\prime}_{e^{i \theta},\lambda}}{F^{\prime}_{e^{i \theta},\lambda}}=\frac{2(1-\beta)\delta(az,\lambda)}{1-\delta(az,\lambda)z}.
	\end{align*}
	A simple computation shows that 
	\begin{align}\label{e-2.7}
		\frac{F^{\prime\prime}_{e^{i \theta},\lambda}}{F^{\prime}_{e^{i \theta},\lambda}}-A(z,\lambda,\beta)&=2(1-\beta)\left(\frac{az(1-\lambda^2)}{(1-\lambda z)(1+\lambda(a-1))z-az^2}\right),
	\end{align}
	\begin{align}\label{e-2.8}
		\frac{F^{\prime\prime}_{e^{i \theta},\lambda}}{F^{\prime}_{e^{i \theta},\lambda}}+B(z,\lambda, \beta)=2(1-\beta)\left(\frac{1-\lambda^2}{(1+\lambda(a-1)z-az^2)(z-\lambda)}\right)
	\end{align}
	and 
	\begin{align}\label{e-2.9}
		\frac{F^{\prime\prime}_{e^{i \theta},\lambda}}{F^{\prime}_{e^{i \theta},\lambda}}-c(z,\lambda, \beta)=\frac{2(1-\beta)(1-\lambda^2)(az-|z|^2)(\lambda a+\bar{z}+\lambda)}{(1-|z|^2)((1+|z|^2)-2\lambda( \real z))(1+\lambda(a-1)z-az^2)}.
	\end{align}
	By replacing $ a $ by $ e^{i\theta} $ ({\it{i.e.}}, substituting $ a=e^{i\theta} $), we obtain
	\begin{align}\label{e-2.10}
		\frac{F^{\prime\prime}_{e^{i\theta},\lambda}}{F^{\prime}_{e^{i\theta},\lambda}}-c(z,\lambda,\beta)=r(z,\lambda,\beta)\frac{e^{i\theta}z\left(1+\lambda\left(e^{-i\theta}-1\right)\bar{z}\right)-a^{-i\theta}z^{-2}}{\left(1+\lambda\left(e^{i\theta}-1\right)z-e^{i\theta}z^2\right)|z|}.
	\end{align}
	From  Lemma \ref{lem-2.11}, we have 
	\begin{align*}
		G^{\prime}(z)=\frac{e^{i\theta}z}{\left(1+\lambda\left(e^{i\theta}-1\right)z-e^{i\theta}z^2\right)^2}.
	\end{align*}
	Thus it is easy to see that
	\begin{align}\label{e-2.12}
		\frac{G^{\prime}(z)}{|G^{\prime}(z)|}=\frac{\left|1+\lambda\left(e^{i\theta}-1\right)-e^{i\theta}z^2\right|^2}{\left(1+\lambda\left(e^{i\theta}-1\right)z-e^{i\theta}z^2\right)^2}\frac{e^{i\theta}z}{|z|}.
	\end{align}
	Applying \eqref{e-2.12} in \eqref{e-2.10}, we obtain
	\begin{align}\label{e-2.13}
		\frac{F^{\prime\prime}_{e^{i\theta},\lambda}}{F^{\prime}_{e^{i\theta},\lambda}}-c(z,\lambda, \beta)=r(z,\lambda,\beta)\frac{\left|1+\lambda\left(e^{i\theta}-1\right)-e^{i\theta}z^2\right|^2}{\left(1+\lambda\left(e^{i\theta}-1\right)z-e^{i\theta}z^2\right)^2}\frac{e^{i\theta}z}{|z|}.
	\end{align}
	Therefore from \eqref{e-2.12} and \eqref{e-2.13} we have 
	\begin{align}\label{e-3.13}
		\frac{F^{\prime\prime}_{e^{i\theta},\lambda}}{F^{\prime}_{e^{i\theta},\lambda}}-c(z,\lambda, \beta)=r(z,\lambda, \beta)\frac{G^{\prime}(z)}{|G^{\prime}(z)|}.
	\end{align}
	Here $ G $ is a starlike function as defined in Lemma \ref{lem-2.11}. 
	Therefore, for any $ z_0\in\mathbb{D}\setminus\{0\} $, the line segment joining $ 0 $ and $ G_0(z_0) $ lies entirely in $ G_0(\mathbb{D}) $. Define $ \gamma_0 $ by
	\begin{align}\label{e-2.15}
		\gamma_0=z(t)=G_0^{-1}\left(tG_0(z_0)\right)\;\; \mbox{for}\,\,0\leq t\leq 1.
	\end{align}
	It is easy to see that 
	\begin{align*}
		G(z(t)) &=2^{-1}e^{i\theta}\left(G_{0}(z(t))\right)^2\\ & =2^{-1}e^{i\theta}\left(tG_{0}(z_0)\right)^2 \\ &=t^2\, 2^{-1}e^{i\theta}G_0(z_0)^2=t^2 \,G(z_0).
	\end{align*}
	Therefore,
	\begin{align}\label{e-2.14}
		G^{\prime}(z(t))z^{\prime}(t)=2tG(z_0)\,\, \mbox{for}\,\, 0\leq t \leq 1.
	\end{align}
	We note that $ z(0)=0 $ and $ z(1)=z_0 $.
	A simple computation using \eqref{e-3.13} shows that
	\begin{align} \label{preeti-e-2.20}
		\log F^{\prime}_{e^{i\theta},\lambda}(z_0)-C(\lambda, \gamma_0, \beta)&=\bigintsss_{0}^{1}\frac{F^{\prime\prime}_{e^{i\theta},\lambda}(z(t))}{F^{\prime}_{e^{i\theta},\lambda}(z(t))}z^{\prime}(t)\, dt-\bigintsss_{0}^{1}c(z(t), \lambda,\beta)z^{\prime}(t)\, dt\\ \nonumber
		&=\bigintsss_{0}^{1}\left(\frac{F^{\prime\prime}_{e^{i\theta},\lambda}(z(t))}{F^{\prime}_{e^{i\theta},\lambda}(z(t))}-c(z(t),\lambda,\beta)\right)z^{\prime}(t)\,dt\\ \nonumber
		&=\bigintsss_{0}^{1}r(z(t),\lambda, \beta)\frac{G^{\prime}(z(t))|z^{\prime}(t)|}{|G^{\prime}(z(t))z^{\prime}(t)|}|z^{\prime}(t)|\, dt\\ \nonumber
		&=\frac{G(z_0)}{|G(z_0)|}\int_{0}^{1}r(z(t),\lambda,\beta)|z^{\prime}(t)|\, dt\\ \nonumber
		&=\frac{G(z_0)}{|G(z_0)|}\, r(\lambda, \gamma_0,\beta).
	\end{align}
	Hence from \eqref{preeti-e-2.20}, we have 
	\begin{align*}
		\left|\log F^{\prime}_{e^{i\theta},\lambda}(z_0)-C(\lambda, \gamma_0,\beta)\right|=R(\lambda,\gamma_0,\beta).
	\end{align*}
	This shows that $ \log F^{\prime}_{e^{i\theta},\lambda}\in\partial \overline{\mathbb{D}}\left(C(\lambda, \gamma_0,\beta), R(\lambda,\gamma_0,\beta)\right) $.
	\\
	
Let $
		V(z_0, \lambda,\beta)=\{\log f^{\prime}(z_0) : f\in \mathcal{F}(\lambda)\}.
	$
	From Corollary \ref{cor-3.6}, it is easy to see that 
	\begin{align*}
		\log F^{\prime}_{e^{i\theta},\lambda}\in\overline{\mathbb{D}}(C(\lambda,\gamma,\beta),R(\lambda,\gamma,\beta))\subset V(z_{0},\lambda,\beta).
	\end{align*}
	Then, we have $ \log F^{\prime}_{e^{i\theta},\lambda}\in\partial V(z_0, \lambda,\beta) $. Suppose there exists $ f $ such that 
	\begin{align*}
		\log F^{\prime}_{e^{i\theta},\lambda}(z_0)=\log f^{\prime}(z_0)
	\end{align*}
	for $z_{0} \in \mathbb{D}$. Then we claim that $ F_{e^{i\theta},\lambda}=f $. We set
	\begin{align*}
		h(t)=\frac{\overline{G(z_0)}}{{G(z_0)}}\left(\frac{f^{\prime\prime}(z(t))}{f^{\prime}(z(t))}-c(z(t), \lambda, \beta)\right)z^{\prime}(t),
	\end{align*}
	where $\gamma_0: z(t)$, $0\leq t \leq 1$ as in \eqref{e-2.15}. Then it is easy to see that $h(t)$ is a continuous function on $[0,1]$ and satisfies
	\begin{align*}
		|h(t)|&=\frac{|\overline{G(z_0)}|}{|{G(z_0)}|}\left|\frac{f^{\prime\prime}(z(t))}{f^{\prime}(z(t))}-c(z(t), \lambda, \beta)\right||z^{\prime}(t)|\\&=\left|\frac{f^{\prime\prime}(z(t))}{f^{\prime}(z(t))}-c(z(t), \lambda, \beta)\right||z^{\prime}(t)|\\&\leq r(z(t), \lambda, \beta)|z^{\prime}(t)| \,\,\,\, \mbox{for} \,\, 0\leq t \leq 1.
	\end{align*}
	Further, we see that
	\begin{align*}
		\int_{0}^{1}{\rm Re} (h(t))dt&=\int_{0}^{1}{\rm Re}\left(\frac{\overline{G(z_0)}}{{G(z_0)}}\left(\frac{f^{\prime\prime}(z(t))}{f^{\prime}(z(t))}-c(z(t), \lambda, \beta)\right)z^{\prime}(t)\right)dt\\&={\rm Re}\left(\frac{\overline{G(z_0)}}{{G(z_0)}}\int_{0}^{1}\frac{f^{\prime\prime}(z(t))}{f^{\prime}(z(t))}z^{\prime}(t)dt-\int_{0}^{1}c(z(t), \lambda, \beta)z^{\prime}(t)dt\right)\\&={\rm Re}\left(\frac{\overline{G(z_0)}}{{G(z_0)}}\left(\log f^{\prime}(z_0)-C(\lambda, \gamma_0,\beta)\right)\right)\\&={\rm Re}\left(\frac{\overline{G(z_0)}}{{G(z_0)}}\left(\log F^{\prime}_{e^{i\theta,\lambda}}(z_0)-C(\lambda, \gamma_0,\beta)\right)\right)\\&=\int_{0}^{1}r(z(t), \lambda, \beta)|z^{\prime}(t)|dt.
	\end{align*}
	Therefore, $ h(t)=r(z(t), \lambda, \beta)|z^{\prime}(t)| $.
	Thus we have 
	\begin{align*}
		\frac{f^{\prime\prime}}{f^{\prime}}=\frac{F^{\prime\prime}_{e^{i\theta,\lambda}}}{F^{\prime}_{e^{i\theta,\lambda}}}
	\end{align*}
	on the curve $ \gamma_0 $. 
	Therefore, we must have 
	$ \log f^{\prime}=\log F^{\prime}_{e^{i\theta,\lambda}} $ which shows that $ f\equiv F_{e^{i\theta,\lambda}} $. 	
This completes the proof.
\end{pf}

\begin{proof}[{\bf Proof of Proposition \ref{pvdev3-th1}}]
We wish to show that the closed curve
$$(-\pi,\pi]\ni \theta \mapsto \log F'_{e^{i\theta},\lambda}(z_0)$$
is a simple curve. Suppose that $\log F'_{e^{i\theta_1,\lambda}}(z_0)=\log F'_{e^{i\theta_2\lambda}}(z_0)$ for some $\theta_1,\theta_2 \in (-\pi,\pi]$ with $\theta_1 \neq \theta_2$. Then in view of \eqref{e-p-5}, we have 
$$
F_{e^{i\theta_1,\lambda}}=F_{e^{i\theta_2\lambda}}.
$$
Using \eqref{main}, it is easy to see that
$$
e^{i\theta_1}z= \tau \left(\dfrac{\omega_{F_{e^{i\theta_1}\lambda}}}{z},\lambda\right)=\tau \left(\dfrac{\omega_{F_{e^{i\theta_2}\lambda}}}{z},\lambda\right)=e^{i\theta_2}z,
$$
which is a contradiction. Thus the curve $(-\pi,\pi]\ni \theta \mapsto \log F'_{e^{i\theta},\lambda}(z_0)$ must be a simple curve.
Since $V(z_0, \lambda,\beta)$ is a compact, convex subset of $\mathbb{C}$ and has non-empty interior, the boundary $\partial V(z_0, \lambda,\beta)$ is a simple closed curve. From Proposition \ref{prop-1}, the curve $\partial V(z_0, \lambda,\beta)$ contains the curve $(-\pi,\pi]\ni \theta \mapsto \log F'_{e^{i\theta},\lambda}(z_0)$. We note that a simple closed curve cannot contain any simple closed curve other than itself. Thus, $\partial V(z_0, \lambda,\beta)$ is given by $(-\pi,\pi]\ni \theta \mapsto \log F'_{e^{i\theta},\lambda}(z_0)$.
\end{proof}

The following figures 1-5 show that the boundary of $V(z_0,\lambda,\beta)$ for the values of $z_0 \in \mathbb{D}\setminus\{0\}, 0 \leq \lambda \leq 1$ and $0 < \beta <1$. We notice that according to Proposition \ref{prop-5.1}, the region bounded by the Jordan curve $\partial V(z_0, \lambda,\beta)$ is compact and convex.
\vspace{2mm}

\begin{figure}[htp]
	\begin{center}
		\includegraphics[width=5.4cm]{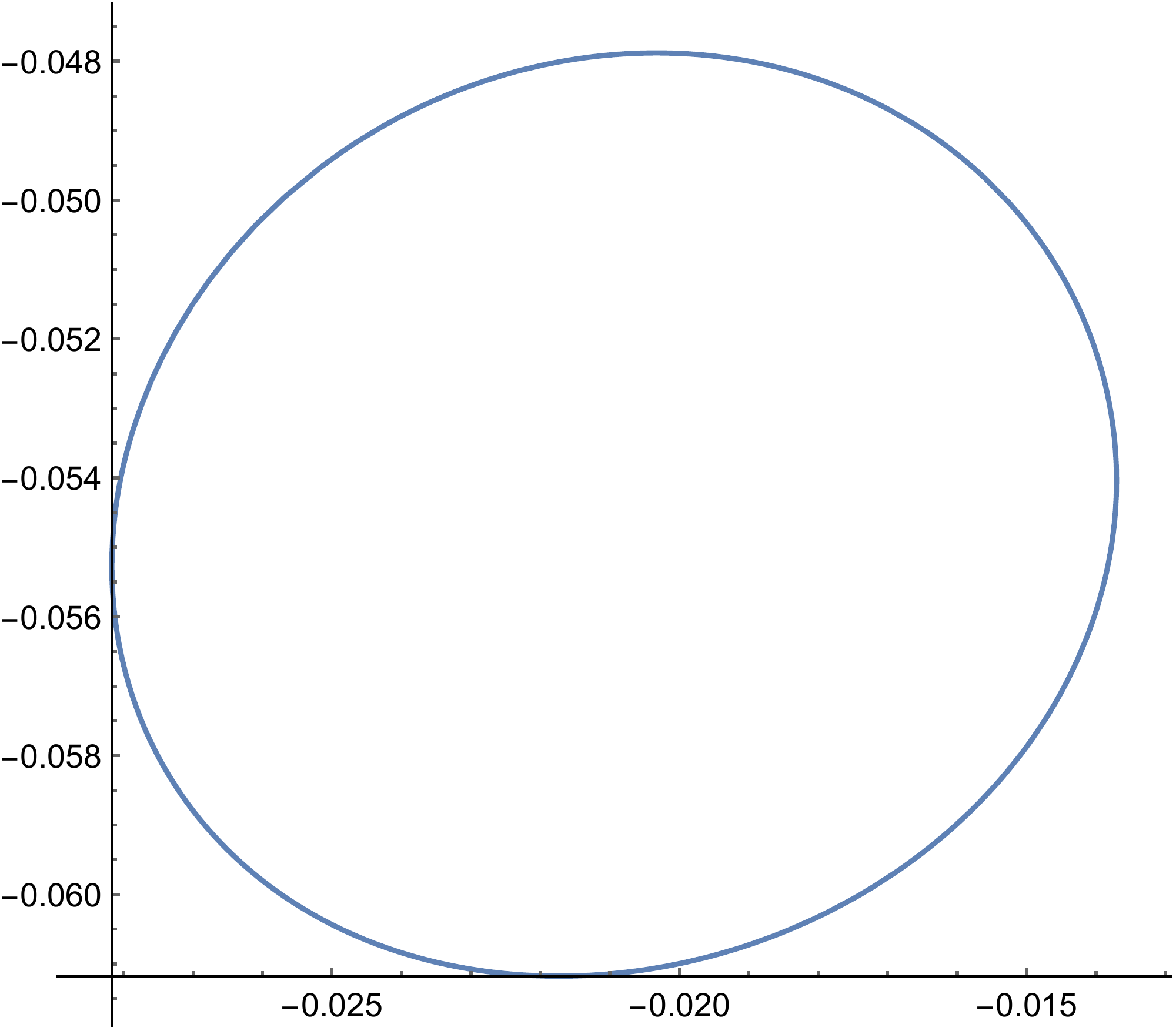}
		\hspace{1cm}
		\includegraphics[width=5.4cm]{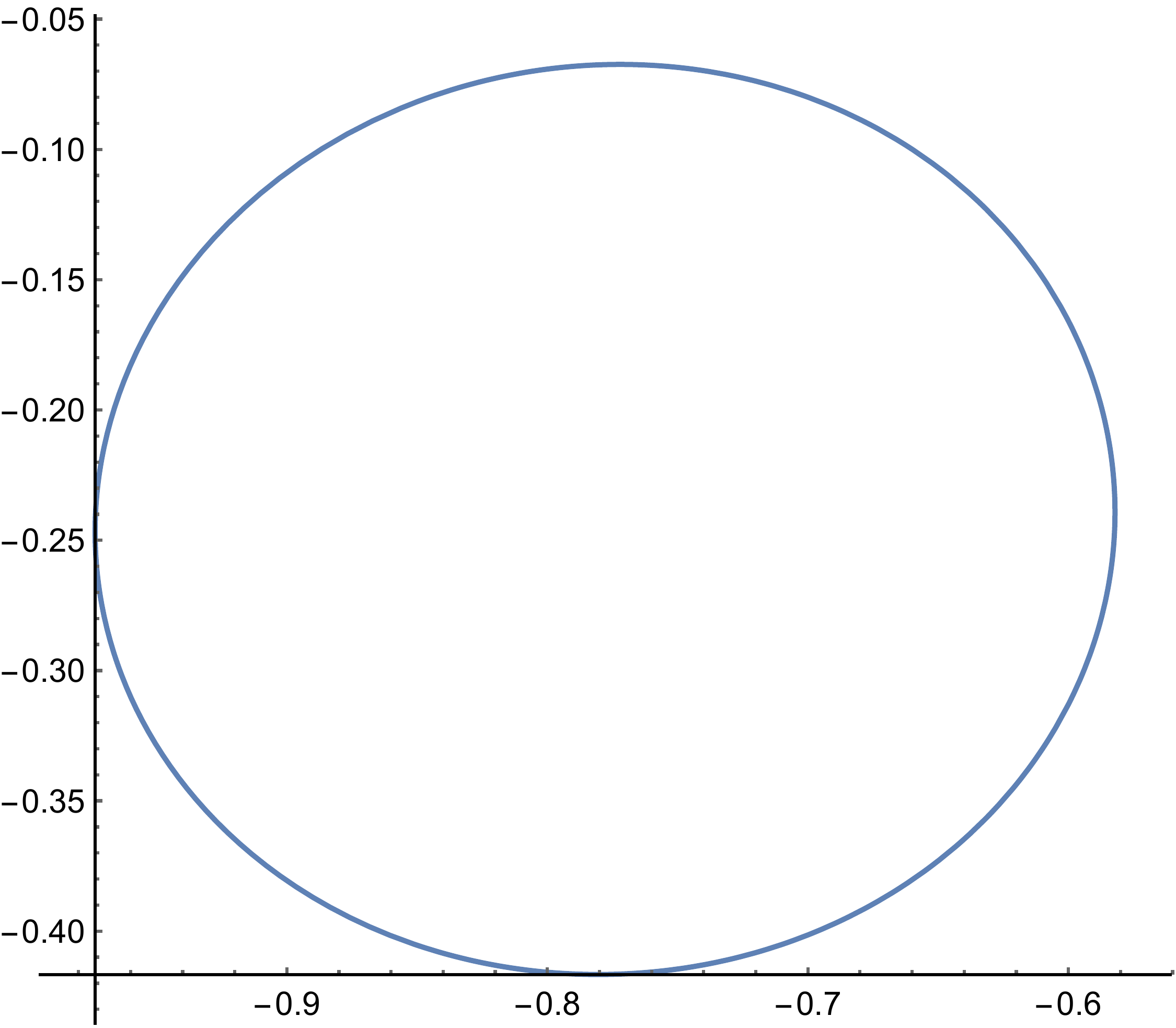}
	\end{center}
	\hspace{1.5cm}$\partial V(z_0,\lambda,\beta)$ \hspace{4.5cm}  $\partial V(z_0,\lambda,\beta)$
	\caption{Region of variability of $\log f'(z_0)$ when $f\in{\mathcal F}(\lambda)$}
\end{figure}
\begin{center}
	$\begin{array}{ll}
		z_0 =0.0778577 + 0.803506i       \hspace{3cm} & z_0 =0.734813 + 0.272699i    \\
		\lambda = 0.886156              & \lambda = 0.751387  \\
		\beta =  0.951601                       & \beta= 0.0264723
	\end{array}$
\end{center}
\vspace{10mm}
\begin{figure}[htp]\label{fig2}
	\begin{center}
		\includegraphics[width=5.4cm]{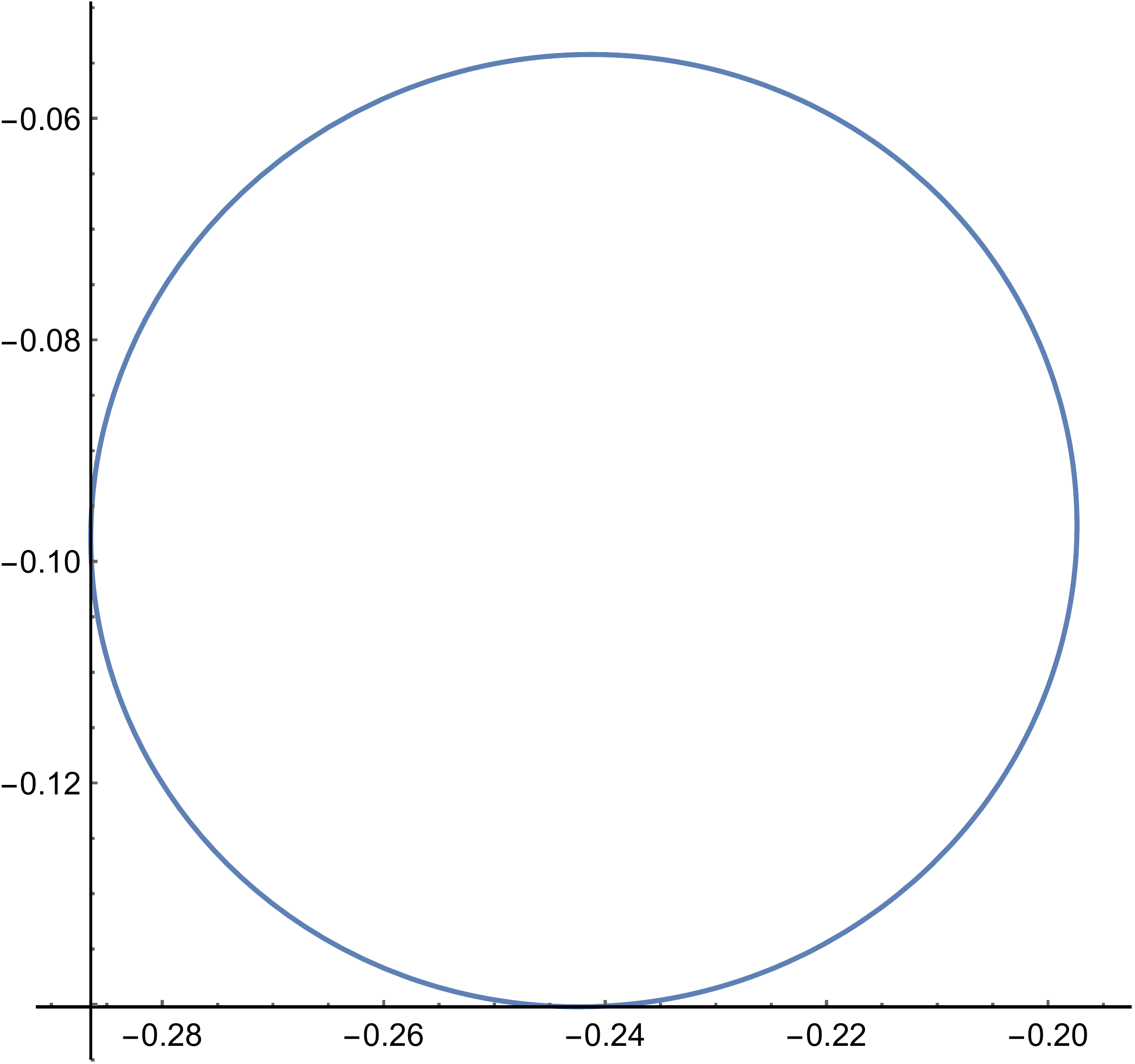}
		\hspace{1cm}
		\includegraphics[width=5.4cm]{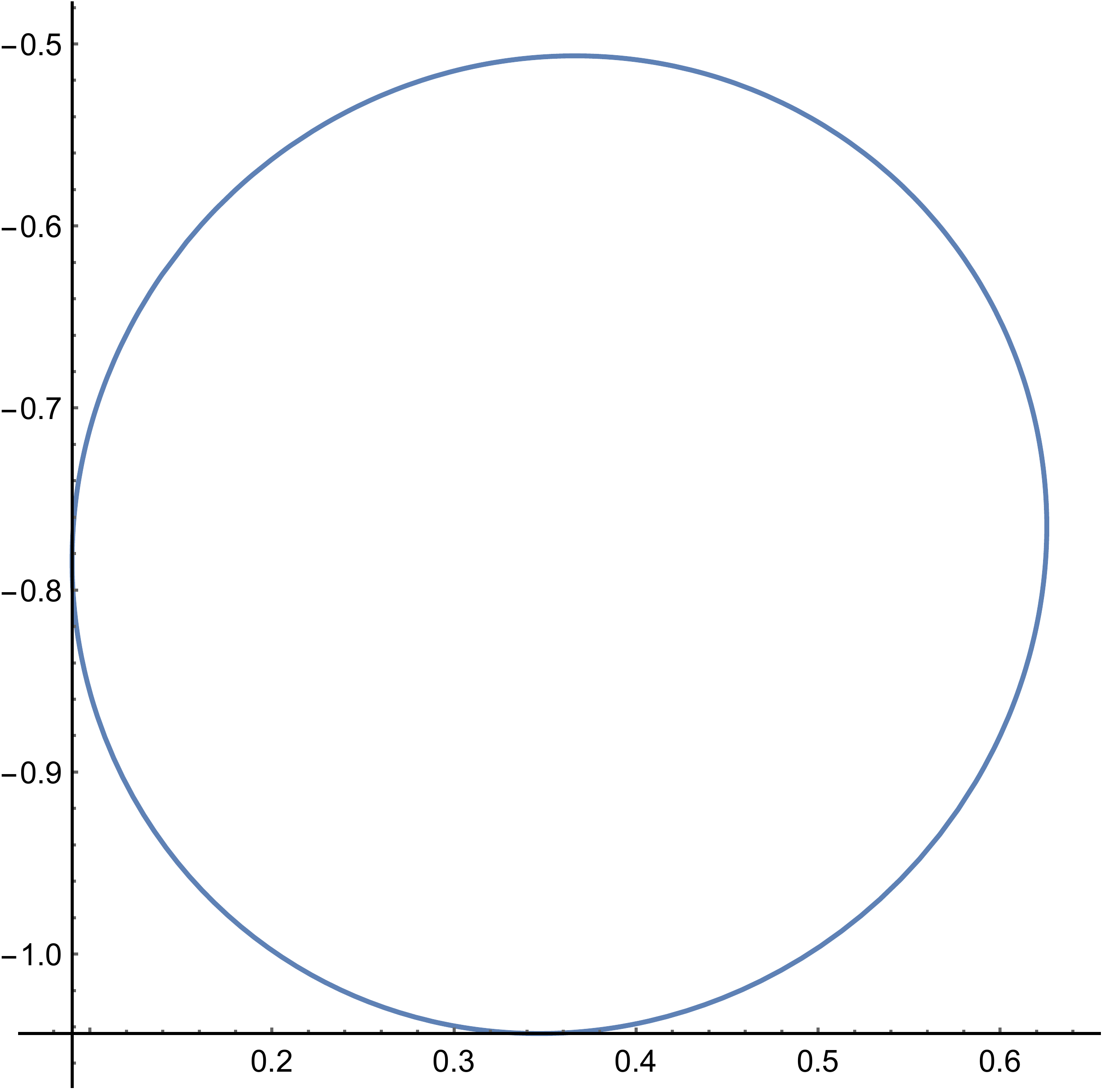}
	\end{center}
	\hspace{1.5cm}$\partial V(z_0,\lambda,\beta)$ \hspace{4.5cm}  $\partial V(z_0,\lambda,\beta)$
	\caption{Region of variability of $\log f'(z_0)$ when $f\in{\mathcal F}(\lambda)$}
\end{figure}
\begin{center}
	$\begin{array}{ll}
		z_0 =0.49963 + 0.23489i       \hspace{3cm} & z_0 =-0.36336 + 0.539691i    \\
		\lambda = 0.716587                & \lambda = 0.71442  \\
		\beta = 0.59288                         & \beta=  0.105526
	\end{array}$
\end{center}
\newpage
\begin{figure}[htp]
	\begin{center}
		\includegraphics[width=5.4cm]{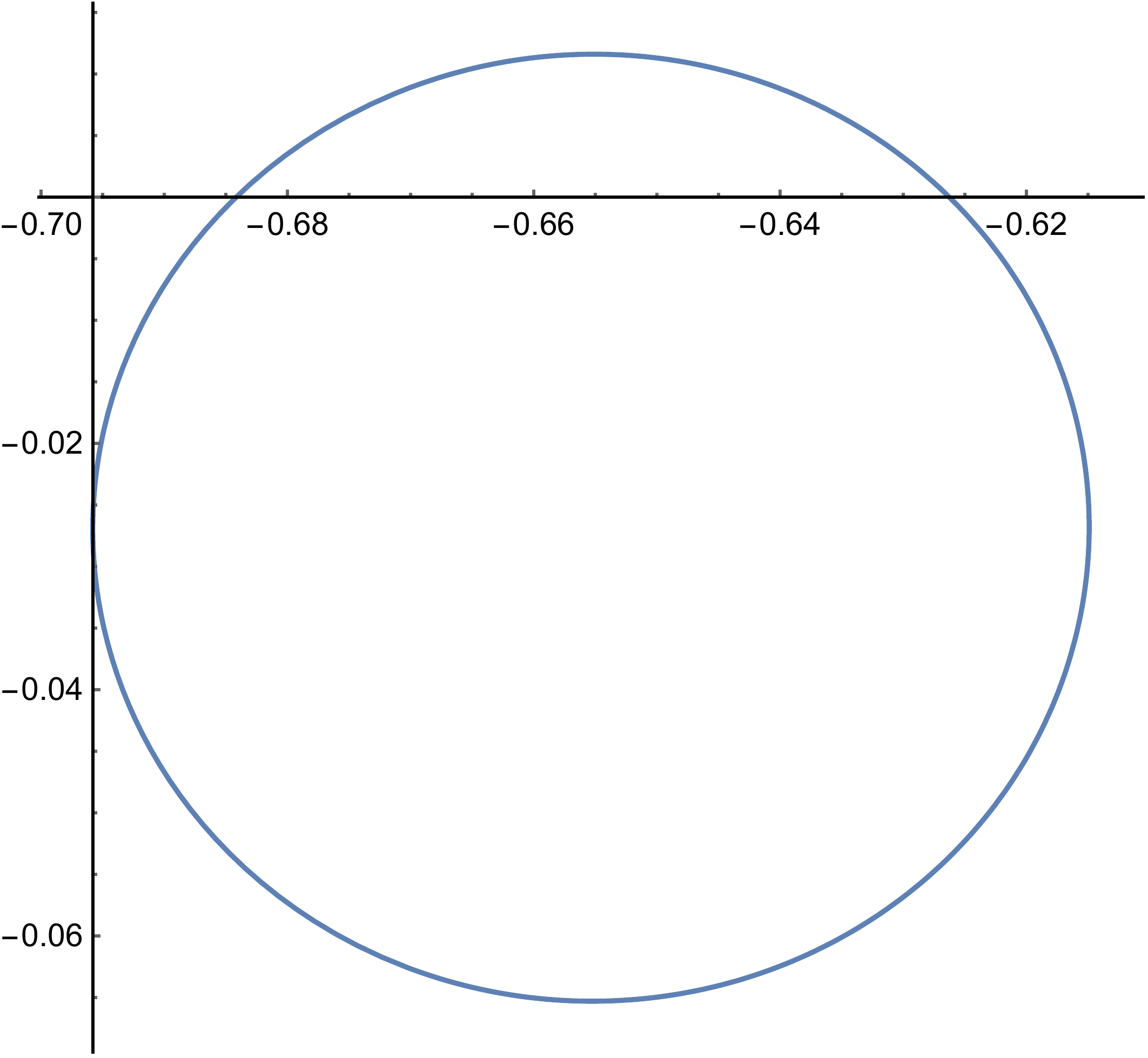}
		\hspace{1cm}
		\includegraphics[width=5.4cm]{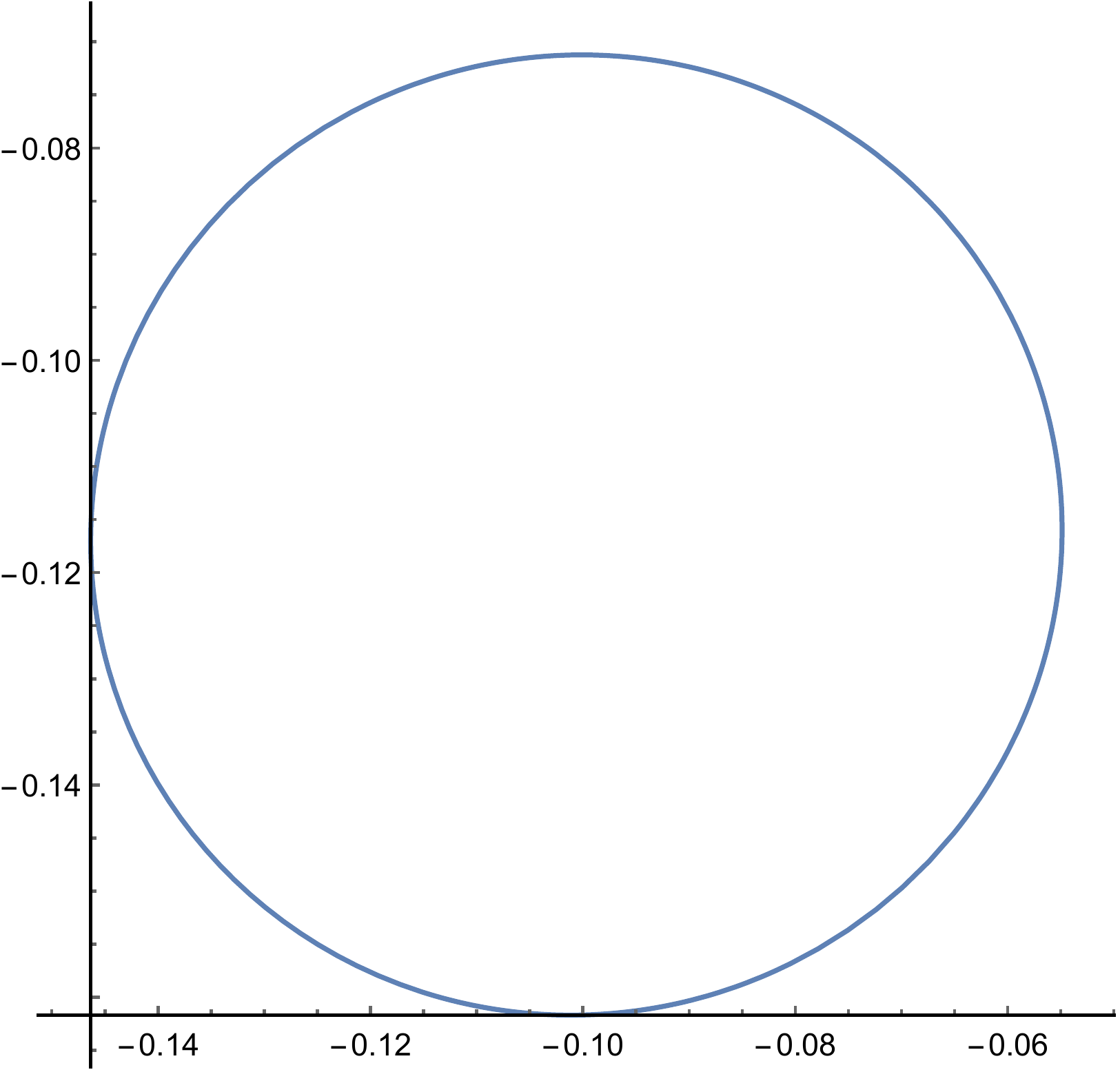}
	\end{center}
	\hspace{1.5cm}$\partial V(z_0,\lambda,\beta)$ \hspace{4.5cm}  $\partial V(z_0,\lambda,\beta)$
	\caption{Region of variability of $\log f'(z_0)$ when $f\in{\mathcal F}(\lambda)$}
\end{figure}
\begin{center}
	$\begin{array}{ll}
		z_0 =0.588369 + 0.0292795i       \hspace{3cm} & z_0 =0.269304 + 0.355388i    \\
		\lambda = 0.884117               & \lambda = 0.529987  \\
		\beta = 0.200402                        & \beta=  0.629449
	\end{array}$
\end{center}
\vspace{10mm}
\begin{figure}[htp]
	\begin{center}
		\includegraphics[width=5.4cm]{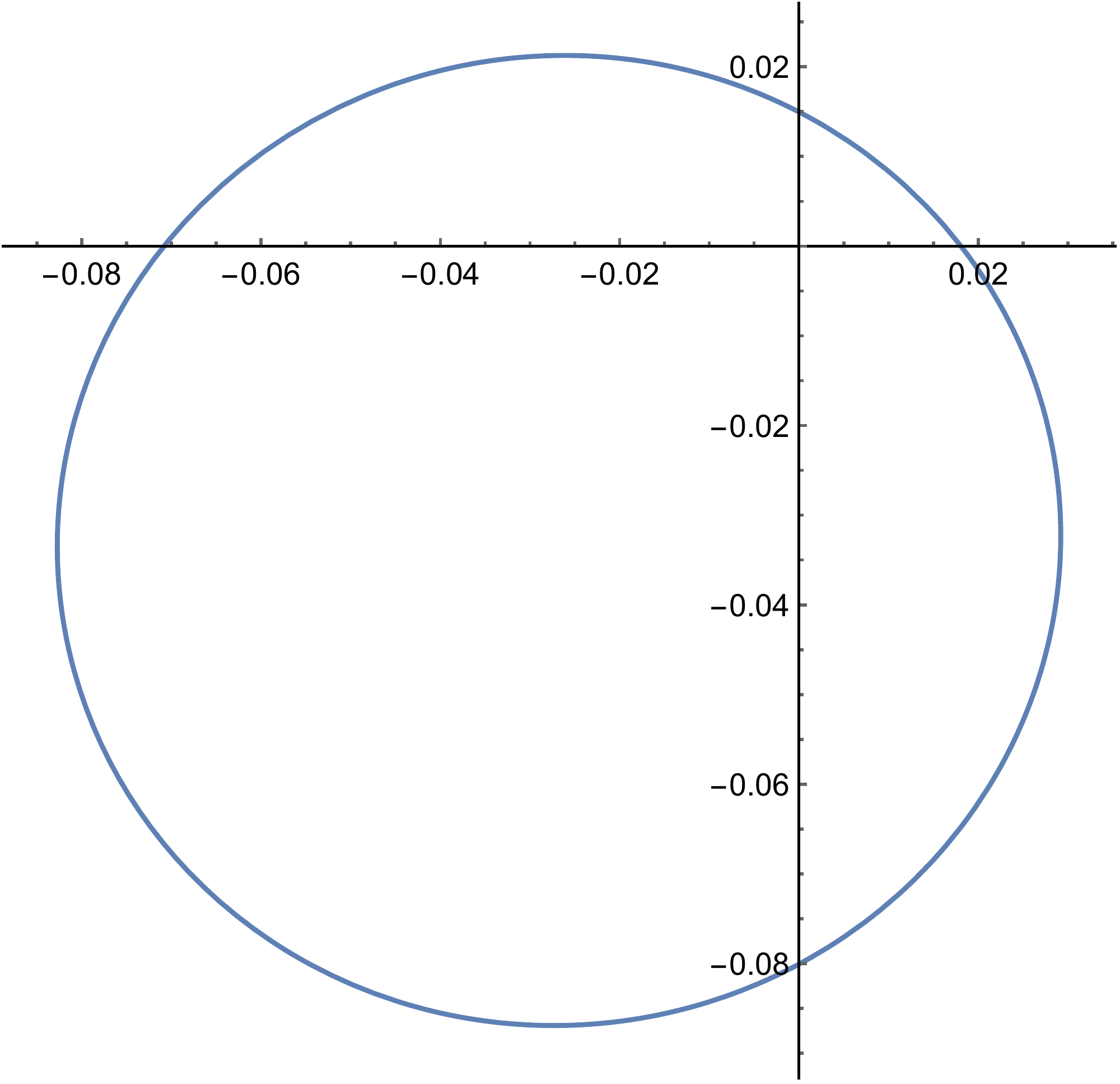}
		\hspace{1cm}
		\includegraphics[width=5.4cm]{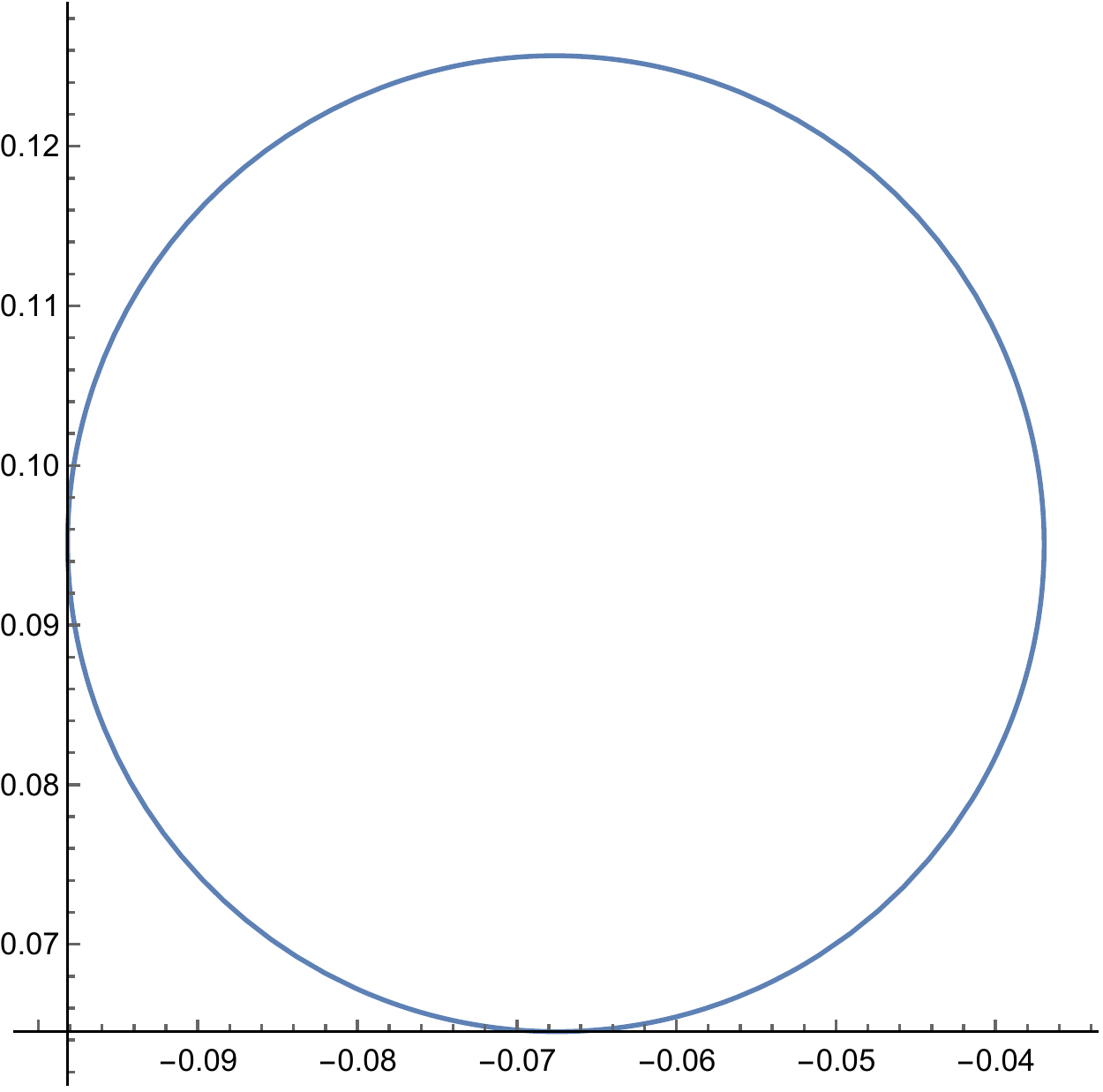}
	\end{center}
	\hspace{1.5cm}$\partial V(z_0,\lambda,\beta))$ \hspace{4.5cm}  $\partial V(z_0,\lambda,\beta)$
	\caption{Region of variability of $\log f'(z_0)$ when $f\in{\mathcal F}(\lambda)$}
\end{figure}
\begin{center}
	$\begin{array}{ll}
		z_0 =0.441545 + 0.437309i       \hspace{3cm} & z_0 =0.138243 - 0.208931i    \\
		\lambda = 0.277434                & \lambda = 0.395356 \\
		\beta = 0.830942                      & \beta=  0.381739
	\end{array}$
\end{center}
\newpage
\begin{figure}[htp]
	\begin{center}
		\includegraphics[width=5.4cm]{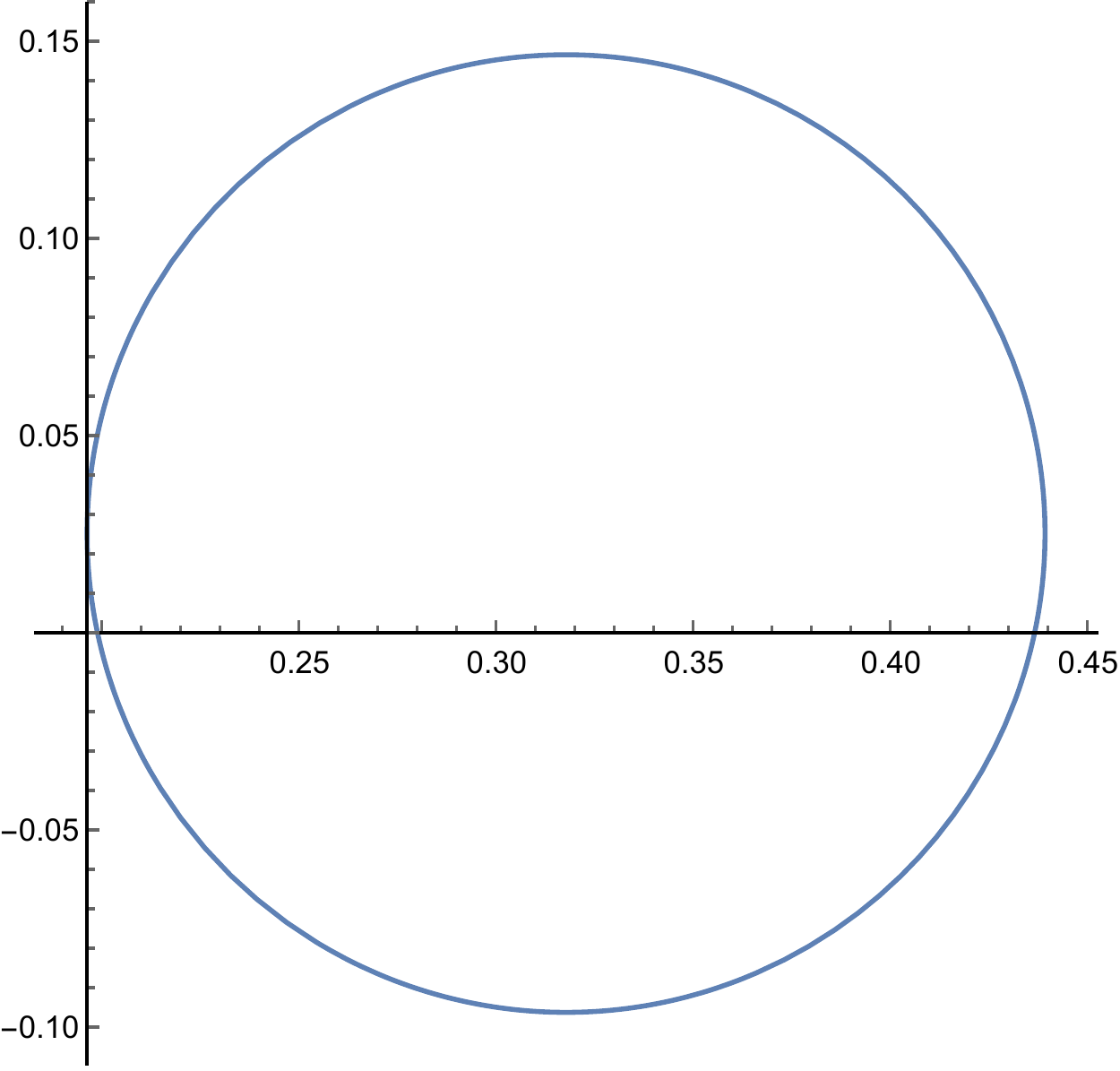}
		\hspace{1cm}
		\includegraphics[width=5.4cm]{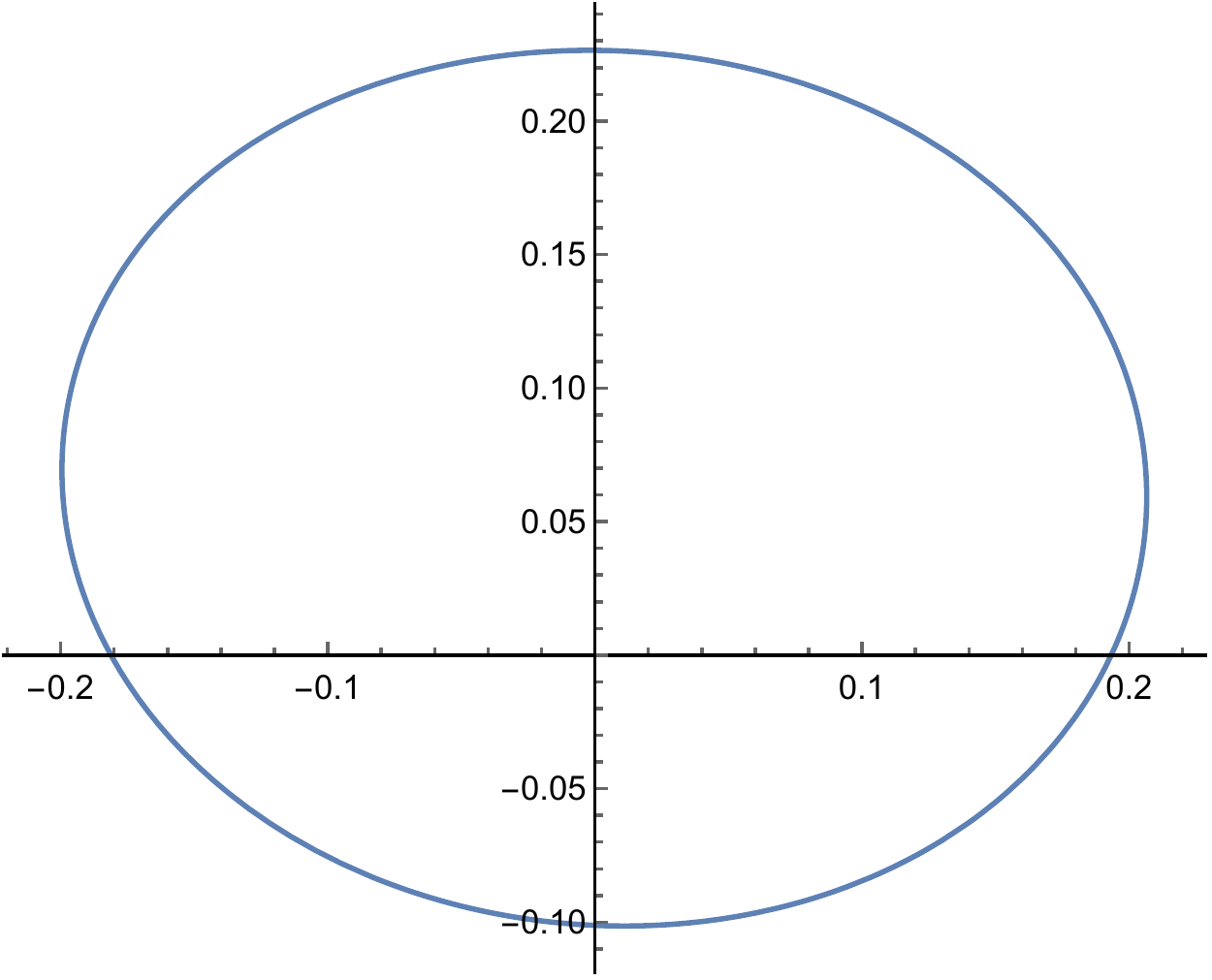}
	\end{center}
	\hspace{1.5cm}$\partial V(z_0,\lambda,\beta)$ \hspace{4.5cm}  $\partial V(z_0,\lambda,\beta)$
	\caption{Region of variability of $\log f'(z_0)$ when $f\in{\mathcal F}(\lambda)$}
\end{figure}
\begin{center}
	$\begin{array}{ll}
		z_0 =-0.533249 - 0.0367948i       \hspace{3cm} & z_0 =0.69143 - 0.596334i    \\
		\lambda = 0.626941              & \lambda = 0.349507  \\
		\beta = 0.593493                    & \beta= 0.762552  
	\end{array}$
\end{center}

\end{document}